\documentclass[a4paper,11pt]{amsart}

\usepackage{amsmath,graphics}
\usepackage{amssymb}
\usepackage{amsfonts}
\usepackage{latexsym}
\usepackage{eucal}
\usepackage{enumerate}
\newtheorem{thm}{Theorem}[section]

\newtheorem{propo}[thm]{Proposition}
\newtheorem{lem}[thm]{Lemma}
\newtheorem{cor}[thm]{Corollary}

\newcommand{\R}{\mathbb{R}}
\newcommand{\C}{\mathbb{C}}
\newcommand{\Z}{\mathbb{Z}}
\newcommand{\N}{\mathbb{N}}

\newcommand{\G}{{\bf G}}

\newcommand{\lt}{{\mathcal L}}

\begin{document}
\bibliographystyle{plain}
\title[Locally constant skew extensions]{Exponentially mixing, locally constant skew extensions of shift maps }
\author{Fr\'ed\'eric Naud}
\address{Laboratoire de Math\'ematiques d'Avignon \\
Campus Jean-Henri Fabre, 301 rue Baruch de Spinoza\\
84916 Avignon Cedex 9. }
\email{frederic.naud@univ-avignon.fr}
 \maketitle
 \begin{abstract}
  We build examples of locally constant $SU_2(\C)$-extensions of the full shift map which are exponentially mixing for the measure of maximal entropy.
 \end{abstract}
 \section{Introduction and result}
 If $T:X\rightarrow X$ is a uniformly hyperbolic map (for example an Anosov diffeomorphism or an expanding map on a smooth manifold), one can produce a {\it skew extension} by considering the new map
 $\widehat{T}:X\times \G \rightarrow X\times \G$ defined by
$$\widehat{T}(x,w):=(Tx,\tau^{-1}(x)w),$$
where $\G$ is a compact connected Lie group and $\tau:X\rightarrow \G$ is a given $\G$-valued function. 
Such a map is an example of partially hyperbolic system.
Given a mixing invariant probability measure $\mu$ for $(X,T)$, the product by the normalized Haar measure $m$ on $\G$ is obviously a $\widehat{T}$-invariant measure for which we can ask natural questions such as mixing, stable ergodicity etc... The qualitative ergodic theory is now fairly well understood, and we refer to the works \cite{Brin1, BurnsWilkinson, Dolgopyat}. The paper of Dolgopyat \cite{Dolgopyat} shows that exponential mixing (for regular enough observables) is generic for extensions of smooth expanding maps, while rapid mixing is also generic for extensions of 
subshifts of finite type. If the map $\tau$ is {\it piecewise constant} on $X$, we call this skew extension {\it locally constant}. It is has been observed first by Ruelle \cite{Ruelle1}, that locally constant extensions of hyperbolic systems cannot be exponentially mixing when $\G$ {\it is a torus}. If $\G=S^1=\R/\Z$, this essentially boils down to the fact that when $\tau$ takes only {\it finitely many values}, Dirichlet box principle shows that one can find sequences of integers 
$N_k$ such that
$$\lim_{k\rightarrow \infty} \sup_{x\in X} \left \vert e^{2i\pi N_k \tau(x)}-1\right \vert=0.$$
For an in-depth study of locally constant {\it toral} extensions, we refer the reader to \cite{Galatolo1} where various rates of polynomial mixing are obtained, depending on the diophantine properties of the values of $\tau$.

In this paper, we will show that surprisingly in the {\it non-commutative case}, it is possible to exhibit a large class of locally constant extensions that are {\it exponentially mixing}. Let us introduce some notations. Let $k\geq 2$ be an integer, and let $\Sigma^+$ be the one-sided shift space
$$\Sigma^+=\{1,\ldots,2k \}^\N,$$
endowed with the shift map $\sigma:\Sigma^+\rightarrow\Sigma^+$ defined by
$$(\sigma \xi)_n=\xi_{n+1},\ \forall\ n\ \in \N.$$
Given $0<\theta<1$, a standard distance $d_\theta$ on $\Sigma^+$ is defined by setting
$$d_\theta(x,\xi)=
\left \{ 1\ \mathrm{if}\ x_0\neq \xi_0  \atop \theta^{N(x,\xi)}\ \mathrm{where}\ N:=\max\{n\geq 1\ :\ x_j=\xi_j\ \forall\ j\leq n \}\  \mathrm{otherwise}. \right.$$
Given a continuous observable $F:\Sigma^+\times \G \rightarrow \C$, which is $d_\theta$-Lipschitz with respect to the first variable, we will use the norm $\Vert F \Vert_{\theta,\G}$ defined by
$$\Vert F \Vert_{\theta,\G}^2:=\sup_{x\in \Sigma^+} \int_{\G} \vert F(x,g)\vert^2dm(g)+\sup_{x\neq \xi}
\frac{1}{d_\theta^2(x,\xi)} \int_{\G} \vert F(x,g)-F(\xi,g)\vert^2dm(g).$$
We will denote by $\mathcal{W}^2_{\theta,\G}$ the completion of the space of continuous functions $F(x,g)$ on 
$\Sigma^+\times \G$, which are $d_\theta$-Lipschitz with respect to $x$ for the $\Vert . \Vert_{\theta,\G}$ norm.

Pick $\tau_1,\ldots,\tau_k \in \G$ and set for all $\xi \in \Sigma^+$,
$$\tau(\xi)=\tau_{\xi_0}\ \mathrm{if}\ \xi_0\in \{1,\ldots,k\},\  \mathrm{and}\ \tau(\xi)=\tau_{\xi_0-k}^{-1}\ \mathrm{if}\ \xi_0\in \{k+1,\ldots,2k\}.$$
Our main result is the following. 
\begin{thm} 
\label{mainth}
Let $\mu$ be the measure of maximal entropy on $\Sigma^+$. Assume that $\G=SU_2(\C)$ and suppose that the group
generated by 
$$\tau_1,\ldots,\tau_k,\tau_1^{-1},\ldots,\tau_k^{-1}$$
is Zariski dense in $\G$ and that $\tau_1,\ldots,\tau_k$ all have algebraic entries. Then there exist $C>0$, $0<\gamma<1$ such that for all $F\in L^2(\Sigma^+\times \G,d\mu dm)$, $G\in \mathcal{W}^2_{\theta,\G}$ we have for all $n\in \N$, 
$$\left \vert \int_{\Sigma^+\times \G} (F\circ \widehat{\sigma}^n) G d\mu dm-\int F d\mu dm \int G d\mu dm \right \vert 
\leq C \gamma^n \Vert F \Vert_{L^2}  \Vert G \Vert_{\theta,\G} .$$
\end{thm}
In other words, we have exponential decay of correlations for H\"older/$L^2$-observables on $\Sigma^+\times \G$, in sharp contrast with the Abelian case. This surprising behaviour will follow from a deep result of Bourgain-Gamburd \cite{BG1}. Extensions to more general groups are possible in view
of the results obtained more recently \cite{BG2, BS}. It is also very likely that the result holds for more general equilibrium measures on $\Sigma^+$,
not just the measure of maximal entropy, and more general subshifts of finite type. This work should be pursued elsewhere.

\bigskip
A corollary of this exponential rate of mixing is the following Central limit theorem for random products in $SU_2(\C)$, which may be of independent interest.
\begin{cor}
 Let $F \in \mathcal{W}^2_{\theta,\G}$ be real valued, and consider the random variable on the probability space $(\Sigma^+\times \G,\mu \times m)$
 $$ Z_n:=\frac{S_n(F)(\xi,g)-n\int\int Fdmd\mu}{\sqrt{n}},$$
 where we have set
 $$S_n(F)(\xi,g)=F(\xi,g)+F(\sigma \xi,\tau^{-1}_{\xi_0}g)+\ldots F( \sigma^{n-1} \xi,\tau^{-1}_{\xi_{n-2}}\ldots \tau^{-1}_{\xi_0}g).$$  Then there exists $\sigma_F\geq 0$ such that as $n\rightarrow \infty$, 
 $Z_n$ converges to a gaussian variable $\mathcal{N}(0,\sigma_F)$ in  law.
 \end{cor}
 The validity of this CLT for generic extensions of hyperbolic systems (including locally constant extensions) and smooth enough $F$ with respect to the $\G$-variable is due to Dolgopyat in \cite{Dolgopyat}. The point of our result is that {\it no smoothness} is required on $F$ with respect to the $\G$-variable, which is unusual in the CLT literature for hyperbolic systems. In particular, if $F(\xi,g)=f(g)$ is in $L^2(\G)$, then the CLT holds.
 This CLT follows directly from absolute summability of correlation functions and the result of Liverani \cite{Liverani}, more details can 
 be found at the end of section $\S$ 2.

\section{Proofs}
\subsection{Bourgain-Gamburd's spectral gap result}
As mentioned earlier, the proof is based on a result of Bourgain-Gamburd \cite{BG1} on the spectral gap of certain Hecke operators acting
on $L^2(\G)$. We start by recalling their result. Given $S=\{\tau_1,\ldots,\tau_k \}\subset \G$, one consider the operator given by
$$T_S(f)(g):=\frac{1}{2k} \sum_{\ell=1}^k \left (f(\tau_\ell g) +f(\tau_\ell^{-1} g)\right).$$
This operator $T_S$ is self-adjoint on $L^2(\G)$ and leaves invariant the one-dimensional space of constant functions. Let $L_0^2(\G)$ be
$$L_0^2(\G):=\left \{f \in L^2(\G)\ :\ \int f dm=0 \right \}.$$
The main result of \cite{BG1} asserts that if the generators $\tau_1,\ldots,\tau_k$ have algebraic entries and the group  
$$\Gamma=\langle \tau_1,\ldots,\tau_k,\tau_1^{-1},\ldots,\tau_k^{-1}\rangle $$
is Zariski dense, then $T_S$ has a spectral gap i.e.
$$\Vert T_S\vert_{L^2_0(\G)} \Vert <1.$$
By self-adjointness, this is equivalent to say that the $L^2$-spectrum of $T_S$ consists of the simple eigenvalue $\{1\}$ while the rest of the spectrum is included in a disc of radius $\rho<1$.
A special case is when $\Gamma$ is a free group, which is enough for many applications. A consequence of the spectral gap property is the following.
If $T_S$ has a spectral gap then there exists $0<\rho<1$ such that for all $n\geq 0$ and $f\in L^2$, 
\begin{equation}
\label{mixing1}
\left \Vert T_S^n(f)-\int fdm\right \Vert_{L^2}\leq 2\rho^n  \Vert f \Vert_{L^2}.
\end{equation}
This deep result of spectral gap is related to the (now solved) Ruziewicz measure problem on invariant means on the sphere. For more details
on the genesis of these problems and explicit examples, we refer the reader to the book of Sarnak \cite{SarnakBook}, chapter 2, and to the paper
\cite{GJS} which was the starting point of \cite{BG1}.
This estimate is the main ingredient of the proof, combined with a decoupling argument and exponential mixing of the measure of maximal entropy.
We recall that given a finite word $\alpha \in \{1,\ldots,2k\}^n$, the {\it cylinder set} $[\alpha]$ of length $\vert \alpha \vert=n$ associated to $\alpha$ is just the set $$[\alpha]:=\{ \xi\in \Sigma^+\ :\ \xi_0=\alpha_1,\ldots,\xi_{n-1}=\alpha_n \}.$$
We will use the following Lemma. 
\begin{lem}
\label{mixing2} For each $\alpha \in \{1,\ldots,2k\}^n$, choose a sequence $\xi_\alpha \in [\alpha]$. let $F$ be an observable with
$$\Vert F\Vert_{\theta,\G}<+\infty.$$ Then we have the bound for all $n\geq 1$,
$$\left \vert     \frac{1}{(2k)^n}\sum_{\vert \alpha \vert=n} \int_G F(\xi_\alpha, g)dm(g)-
\int \int F d\mu dm \right \vert \leq  \Vert F\Vert_{\theta,\G} \theta^n.$$
\end{lem}
\noindent {\it Proof.} This is a rephrasing of the fact that the measure of maximal entropy is exponentially mixing for H\"older observables on $\Sigma^+$.
Indeed, recall that in our case, the measure of maximal entropy $\mu$ is just the Bernoulli measure such that
$$\mu([\alpha])=\frac{1}{(2k)^n}$$
when $\vert \alpha\vert=n$. By Schwarz inequality, for all $x\in [\alpha]$, we have
$$\left \vert \int_\G F(x,g)dm(g)-\int_\G F(\xi_\alpha,g)dm(g) \right\vert \leq \Vert F\Vert_{\theta,\G}d_\theta(\xi_\alpha,x)$$
$$\leq \Vert F\Vert_{\theta,\G} \theta^n.$$
Writing
$$ \int \int F d\mu dm=\sum_{\vert \alpha \vert=n} \int_{[ \alpha ]} \int_{\G}F(x,g)dm(g),$$
we deduce
$$\left \vert  \sum_{\vert \alpha \vert=n} \mu([\alpha])\int_G F(\xi_\alpha, g)dm(g)-
\int \int F d\mu dm \right \vert \leq  \Vert F\Vert_{\theta,\G} \theta^n,$$
and the proof is done. $\square$

\bigskip
In the sequel, we will use the following notation: given a finite word $\alpha \in \{1,\ldots,2k\}^n$ and $\xi \in \Sigma^+$, we will denote by
$\alpha\xi$ the concatenation of the two words i.e. the new sequence $\alpha\xi \in \Sigma^+$ such that
$(\alpha \xi)_j=\alpha_{j+1}$ for $j=0,\ldots,n-1$ and $\sigma^n(\alpha \xi)=\xi$. We recall that given $f,g\in C^0(\Sigma^+)$, we have the {\it transfer operator identity}
$$\int_{\Sigma^+}(f\circ \sigma^n) g d\mu=\int_{\Sigma^+} f \lt^n(g) d\mu,$$
where we have
$$\lt^n(g)(\xi)=\frac{1}{(2k)^n}\sum_{\vert \alpha \vert =n} g(\alpha \xi).$$
This identity follows straightforwardly from the $\sigma$-invariance of the measure $\mu$ and its value on cylinder sets. Notice that the operator $\lt$ is {\it normalized} i.e. satisfies
$$\lt( {\bf 1})={\bf 1},$$
which will be used throughout all the computations below. Notice also that using the above notations,
we have for all $f \in L^2(\G)$, for all $x\in \Sigma^+$,
$$T_S(f)(g)=\frac{1}{2k}\sum_{\ell=1}^{2k}f(\tau(\ell x)g),$$
while
$$T^n_S(f)(g)=\frac{1}{(2k)^n}\sum_{\vert \beta \vert=n} f(\tau(\beta x)\ldots \tau(\beta_{n-1}\beta_{n}x)\tau(\beta_{n}x)g).$$
The fact that $\tau(x)$ depends only on the first coordinate of $x$ is critical in the above identities.

\subsection{Main proof}
We now move on to the proof of the main result. Let $F,G \in C^0(\Sigma^+\times \G)$, and compute the correlation function:
$$\int\int (F\circ \widehat \sigma^n) G d\mu dm$$
$$=\int \int F(\sigma^n x, \tau^{-1}(\sigma^{n-1}x)\ldots\tau^{-1}(\sigma x) \tau^{-1}(x)g)G(x,g)d\mu(x)dm(g).$$
By using Fubini and translation invariance of the Haar measure we get
$$\int\int (F\circ \widehat \sigma^n) G d\mu dm=\int\int F(\sigma^n x,g)G(x,\tau^{(n)}(x)g)dm(x)dm(g),$$
where
$$\tau^{(n)}(x)=\tau(x)\tau(\sigma x)\ldots \tau(\sigma^{n-1} x).$$ Again by Fubini and the transfer operator formula we get
$$ \int\int (F\circ \widehat \sigma^n) G d\mu dm=\int\int F \widehat{\lt}^n(G)d\mu dm,$$
where
$$ \widehat{\lt}^n(G)(x,g)=\frac{1}{(2k)^n}\sum_{\vert \alpha \vert=n}G(\alpha x, \tau^{(n)}(\alpha x)g).$$
The main result will follow from the following estimate.
\begin{propo}
\label{mixing3}
There exist $C>0$ and $0<\gamma<1$ such that for all $n\geq 1$,
$$\sup_{x\in \Sigma^+} \left \Vert  \widehat{\lt}^n(G)(x,g)-\int \int Gd\mu dm \right \Vert_{L^2(\G)}\leq C \Vert G \Vert_{\theta, \G} \gamma^n.$$
\end{propo}
\noindent Indeed, write 
$$\int\int F \widehat{\lt}^n(G)d\mu dm-\int\int Fd\mu dm \int\int Gd\mu dm$$
$$= \int \int F\left(\widehat{\lt}^n(G)-\int\int G d\mu dm\right)d\mu dm,$$
and use Schwarz inequality combined with the above estimate to get the conclusion of the main theorem. Let us prove Proposition \ref{mixing3}.
Writing $n=n_1+n_2$, we get 
$$\widehat{\lt}^n(G)(x,g)=\frac{1}{(2k)^{n_1}} \sum_{\vert \alpha \vert=n_1} \frac{1}{(2k)^{n_2}} \sum_{\vert \beta \vert=n_2} 
G(\alpha\beta x,\tau^{(n)}(\alpha \beta x)g).$$
Observe that we have (since $\tau$ depends only on the first coordinate)
$$\tau^{(n)}(\alpha \beta x)=\underbrace{\tau(\alpha\beta x)\ldots \tau(\alpha_{n_1}\beta x)}_{depends\  only\ on\  \alpha}
\underbrace{\tau(\beta x)\ldots \tau(\beta_{n_2} x)}_{depends\  only\ on\  \beta}.$$
For all word $\alpha \in \{1,\ldots,2k\}^{n_1}$, we choose $\xi_\alpha \in [\alpha]$, and set 
$$G_\alpha(g):=G(\xi_\alpha,\tau_{\alpha_1}\ldots \tau_{\alpha_{n_1}}g).$$
We now have
$$\widehat{\lt}^n(G)(x,g)=\frac{1}{(2k)^{n_1}} \sum_{\vert \alpha \vert=n_1}T_S^{n_2}(G_\alpha)(g) 
+\mathcal{R}_n(x,g),$$
where the "remainder" $\mathcal{R}_n(x,g)$ is 
$$\mathcal{R}_n(x,g)=$$
$$\frac{1}{(2k)^{n_1}} \sum_{\vert \alpha \vert=n_1} \frac{1}{(2k)^{n_2}} \sum_{\vert \beta \vert=n_2} \left ( G(\alpha\beta x,\tau^{(n)}(\alpha \beta x)g)
-G_\alpha( \tau_{\beta_1}\ldots \tau_{\beta_{n_2}}g)\right).$$
Note that for all $x\in \Sigma^+$, we have by translation invariance of Haar measure, 
$$\left \Vert G(\alpha\beta x,\tau^{(n)}(\alpha \beta x)g)
-G_\alpha(\tau_{\beta_1}\ldots \tau_{\beta_{n_2}} g) \right \Vert_{L^2(\G)}\leq \Vert G \Vert_{\theta, \G} d_\theta(\alpha\beta x,\xi_\alpha),$$
and thus
$$\Vert \mathcal{R}_n(x,g)\Vert_{L^2(\G)}\leq \theta^{n_1}\Vert G\Vert_{\theta,\G}.$$ On the other hand, using Lemma \ref{mixing2},
we have 
$$\left \Vert \frac{1}{(2k)^{n_1}} \sum_{\vert \alpha \vert=n_1}T_S^{n_2}(G_\alpha)(g) -\int\int Gd\mu dm \right\Vert_{L^2(\G)}$$
$$\leq  \frac{1}{(2k)^{n_1}} \sum_{\vert \alpha \vert=n_1} \left \Vert T_S^{n_2}(G_\alpha)-\int_{\G} G(\xi_\alpha,g)dm(g) \right \Vert_{L^2(\G)}
+O(\Vert G \Vert_{\theta,\G}\theta^{n_1}).$$
By the spectral gap property (\ref{mixing1}),
$$\left \Vert T_S^{n_2}(G_\alpha)-\int_{\G} G(\xi_\alpha,g)dm(g) \right \Vert_{L^2(\G)}\leq 2\rho^{n_2} \Vert G_\alpha \Vert_{L^2(\G)}\leq 
2\rho^{n_2} \Vert G \Vert_{\theta,\G},$$
therefore we have obtained, uniformly in $x\in \Sigma^+$, 
$$ \left \Vert  \widehat{\lt}^n(G)(x,g)-\int \int Gd\mu dm \right \Vert_{L^2(\G)}=O\left( \Vert G\Vert_{\theta,G}(\theta^{n_1}+\rho^{n_2}) \right), $$
and the proof ends by choosing $n_1=[n/2]$, $n_2=n-n_1$. $\square$

\bigskip
In the proof, we have used no specific information about the group $\G$, except that $T_S$ has a spectral gap. Therefore the main theorem extends without modification to the case of $\G=SU_d(\C)$ and more generally any compact connected simple Lie group by \cite{BS}. 

We also point out that the rate of mixing obtained is $O\left ( (\max\{\sqrt{\theta}, \sqrt{\rho} \})^n \right)$, which can be made explicit if $T_S$ has an explicit
spectral gap, see \cite{GJS} for some examples arising from quaternionic lattices. On the other hand, if
$$\Gamma=\langle \tau_1,\ldots,\tau_k,\tau_1^{-1},\ldots,\tau_k^{-1}\rangle $$
is free, it follows from a result of Kesten \cite{Kesten} that the $L^2$-spectrum of $T_S$ contains the full continuous segment 
$$\left[-\frac{\sqrt{2k-1}}{k},+\frac{\sqrt{2k-1}}{k}\right],$$
 which suggests that the rate of mixing cannot exceed $O\left( \left(\frac{\sqrt{2k-1}}{k}\right)^n \right)$.
\subsection{Central Limit Theorem} The CLT, has stated in the introduction follows from the paper of Liverani \cite{Liverani}, section $\S 2$ on non invertible,
onto maps. Following his notations we have for all $F\in L^2(dmd\mu)$,
$$T(F)=F\circ \widehat{\sigma},$$
while the $L^2$-adjoint $T^*$ is $T^*=\widehat{\lt}$.
Given $F\in \mathcal{W}^2_{\theta,\G}$ with $\int \int F dmd\mu=0$, we know that :
\begin{enumerate}
 \item we have by Theorem \ref{mainth} $$ \sum_{\ell=0}^\infty \left \vert \int \int (F\circ \widehat{\sigma}^\ell) F   d\mu dm \right \vert <\infty,$$
 \item by Proposition \ref{mixing3}, the series
 $$\sum_{\ell \in \N} \widehat{\lt}^\ell (F),$$
converges absolutely almost surely (actually in $L^1(dmd\mu)$).  
\end{enumerate}
We can therefore apply Therorem 1.1 from \cite{Liverani}, which says that the CLT holds and that the
variance $\sigma_F$ is vanishing if and only of $F$ is a coboundary. Alternatively, since we have a spectral gap for $\widehat{\lt}$
\footnote{Proposition \ref{mixing2} can indeed be refined to show that $\widehat{\lt}:\mathcal{W}^2_{\theta,\G}\rightarrow \mathcal{W}^2_{\theta,\G}$ has a spectral gap, but we have chosen here the shortest path.} one could also use a spectral method and perturbation theory to prove the CLT, see \cite{Gouezel} for a survey on this approach for proving central limit theorems.

\bigskip
\noindent{\bf Acknowledgments.} Fr\'ed\'eric Naud is supported by ANR GeRaSic and Institut universitaire de France.


\begin{thebibliography}{10}

\bibitem{BS}
Yves Benoist and Nicolas de~Saxc\'e.
\newblock A spectral gap theorem in simple {L}ie groups.
\newblock {\em Invent. Math.}, 205(2):337--361, 2016.

\bibitem{BG2}
J.~Bourgain and A.~Gamburd.
\newblock A spectral gap theorem in {${\rm SU}(d)$}.
\newblock {\em J. Eur. Math. Soc. (JEMS)}, 14(5):1455--1511, 2012.

\bibitem{BG1}
Jean Bourgain and Alex Gamburd.
\newblock On the spectral gap for finitely-generated subgroups of {$\rm
  SU(2)$}.
\newblock {\em Invent. Math.}, 171(1):83--121, 2008.

\bibitem{Brin1}
M.~I. Brin.
\newblock The topology of group extensions of {$C$}-systems.
\newblock {\em Mat. Zametki}, 18(3):453--465, 1975.

\bibitem{BurnsWilkinson}
Keith Burns and Amie Wilkinson.
\newblock Stable ergodicity of skew products.
\newblock {\em Ann. Sci. \'Ecole Norm. Sup. (4)}, 32(6):859--889, 1999.

\bibitem{Dolgopyat}
Dmitry Dolgopyat.
\newblock On mixing properties of compact group extensions of hyperbolic
  systems.
\newblock {\em Israel J. Math.}, 130:157--205, 2002.

\bibitem{Galatolo1}
Stefano Galatolo, J\'er\^ome Rousseau, and Benoit Saussol.
\newblock Skew products, quantitative recurrence, shrinking targets and decay
  of correlations.
\newblock {\em Ergodic Theory Dynam. Systems}, 35(6):1814--1845, 2015.

\bibitem{GJS}
Alex Gamburd, Dmitry Jakobson, and Peter Sarnak.
\newblock Spectra of elements in the group ring of {${\rm SU}(2)$}.
\newblock {\em J. Eur. Math. Soc. (JEMS)}, 1(1):51--85, 1999.

\bibitem{Gouezel}
S\'ebastien Gou\"ezel.
\newblock Limit theorems in dynamical systems using the spectral method.
\newblock In {\em Hyperbolic dynamics, fluctuations and large deviations},
  volume~89 of {\em Proc. Sympos. Pure Math.}, pages 161--193. Amer. Math.
  Soc., Providence, RI, 2015.

\bibitem{Kesten}
Harry Kesten.
\newblock Symmetric random walks on groups.
\newblock {\em Trans. Amer. Math. Soc.}, 92:336--354, 1959.

\bibitem{Liverani}
Carlangelo Liverani.
\newblock Central limit theorem for deterministic systems.
\newblock In {\em International {C}onference on {D}ynamical {S}ystems
  ({M}ontevideo, 1995)}, volume 362 of {\em Pitman Res. Notes Math. Ser.},
  pages 56--75. Longman, Harlow, 1996.

\bibitem{Ruelle1}
David Ruelle.
\newblock Flots qui ne m\'elangent pas exponentiellement.
\newblock {\em C. R. Acad. Sci. Paris S\'er. I Math.}, 296(4):191--193, 1983.

\bibitem{SarnakBook}
Peter Sarnak.
\newblock {\em Some applications of modular forms}, volume~99 of {\em Cambridge
  Tracts in Mathematics}.
\newblock Cambridge University Press, Cambridge, 1990.

\end{thebibliography}
\end{document}